\documentstyle[12pt]{article}
\def\v{\mbox{\rm Vol}}
\begin{document}
\begin{center}
{\Large \bf On the volumes of complex hyperbolic manifolds }
\\ {\Large \bf with cusps}

\bigskip
 {\large \bf Jun-Muk Hwang}
\footnote{Supported  by the Korea Research Foundation Grant
(KRF-2002-070-C00003).}
 \end{center}

\bigskip
\begin{abstract}
We study the problem of bounding the number of cusps of a complex
hyperbolic manifold in terms of its volume. Applying
algebro-geometric methods  using Mumford's work on toroidal
compactifications and its generalization due to N. Mok and W.-K.
To, we get a bound which is considerably better than those
obtained previously by methods of geometric topology.
\end{abstract}

\bigskip
{\bf MSC Number:} 32Q45

\section{Introduction}
There have been some interests on the problem of bounding the
number of cusps of a complex hyperbolic manifold in terms of its
volume, as a generalization of the corresponding problem for a
real hyperbolic manifold.  We refer the readers to [3], [5] and
the references therein for the historical background  and the
motivation for studying problems of this type from the view-point
of geometric topology.

It seems that the following bounds of John R. Parker's are the
best published results on this problem.

\medskip
{\bf Theorem 1 [5, Theorem D and Theorem F]} {\it Let $X$ be an
$n$-dimensional complex hyperbolic manifold of finite volume. Let
$k$ be the number of cusps of $X$ and let $\v(X)$ be the volume of
$X$ with respect to the Bergmann metric with holomorphic sectional
curvature $-1$. Then $$ \frac{\v(X)}{k} \geq \frac{2^{n-1}}{n (6
\pi)^{2n^2 - 3n +1}}.$$ When $n=2$,
$$\frac{\v(X)}{k} \geq \frac{2}{3}.$$}

\medskip
 The method used in [5], based on the earlier work of [3], is
motivated by the corresponding method in the study of real
hyperbolic manifolds. More precisely, these authors constructed
certain disjoint neighborhoods of the cusps whose volumes can be
estimated.

\medskip
The goal of this paper is to explain  a completely different
approach to the problem, using techniques of  algebraic geometry.
To state our result, let $$P(\ell) := \frac{(n \ell + n +
\ell)!}{n! (n \ell + \ell)!}.$$

\medskip
{\bf Theorem 2} {\it  In the notation of Theorem 1, for $n \geq
2$,
$$\frac{\v(X)}{k} \geq
\frac{(4 \pi)^n}{ n! (P(4) - P(2))} (1- \frac{n+1}{P(4) -
P(2)}).$$}

\medskip
Note that the right hand side is at least $
 \frac{2^{2n-1} \pi^n}{(5n+4)^n}$ which is considerably better
 than Theorem 1.  For $n=2$, $P(4) - P(2) = 63$ and
  the right hand side
 is $$ \frac{(4 \pi)^2}{2 \cdot 63}(1- \frac{3}{63}) =
 \frac{160}{1323} \pi^2 \geq 1.19...,$$
which is better than Theorem 1. Note that our argument is uniform
in all dimensions $\geq 2$, while the case $n=2$ in Parker's work
was obtained by a special argument which did not apply in higher
dimensions.

\medskip
Theorem 2 is  obtained by examining the dimensions of the spaces
of certain cusp forms. The proof depends essentially on the
existence of a toroidal compactification of $X$ and its metric
property which was established by Mumford [4] for $X$ defined by
an arithmetic group and generalized to arbitrary $X$ by N. Mok and
W.-K. To [7]. Excepting these results, we only need standard
methods of algebraic geometry.

\medskip
Yum-Tong Siu told us that one may be able to get a bound of the
above type also by the differential geometric method used in [6].
It is not clear however whether the resulting bound would be as
good as ours.

\section{Results from  toroidal compactifications}

In this section, we will recall some basic facts about toroidal
compactifications which we need for the proof of Theorem 2.

\medskip
Throughout, $X$ denotes a complex hyperbolic manifold of dimension
$n \geq 2$ with finite volume.
 Denote by $X^*$ the minimal compactifcation of $X$, which was constructed
by Baily-Borel [2] for $X$ defined by an arithmetic group and by
Siu-Yau [6] for arbitrary $X$. The complement $X^*\setminus X$
consists of $k$ cusp points, which we denote by $$ X^* \setminus X
= \{ Q_1, \ldots, Q_k\}.$$
  $X^*$ is a normal projective
variety and there exists an ample line bundle $K_{X^*}$ extending
the canonical bundle of $X$.

 Denote by $\bar{X}$  a toroidal compactification of $X$, which
 was constructed by Mumford et al. [1] for $X$
defined by an arithmetic group and by Mok for arbitrary $X$ as
explained in [7, p.61].
  $\bar{X}$ is a smooth projective
 variety and the
 complement $ \bar{X} \setminus X$ is a smooth divisor $E$ with $k$ components, which we denote by
$$ \bar{X} \setminus X = E = E_1 \cup \cdots \cup E_k.$$ Each component $E_i$
is an abelian variety of dimension $n-1$ whose normal bundle in
$\bar{X}$ is a negative line bundle, as described in [7,
pp.61-62]. There is a canonical morphism $$\psi: \bar{X}
\rightarrow X^*$$ which contracts each $E_i$ to a cusp point
$Q_i$. Let us denote by $L$ the nef and big line bundle $\psi^*
K_{X^*}.$ Then by [4, Proposition 3.4 (b)], $$L  = K_{\bar{X}} +
E.$$

The key property of $L$ is that the Bergman metric on $X$ induces
a singular metric on $L$ which is good in the sense of [4, Section
1]. This was proved by [4, Main Theorem 3.1 and Proposition 3.4
(b)] for $X$ defined by an arithmetic group and by [7, Section 2]
for arbitrary $X$.  This implies  Hirzebruch proportionality [4,
Theorem 3.2]. One special case we need is the following.

\medskip {\bf Proposition 1 [4, Theorem 3.2]}
$$ \v(X) = \frac{(4 \pi)^n}{n! (n+1)^{n}} L^n .$$

\medskip
This is not exactly [4, Theorem 3.2] because Mumford uses
different normalization of the metric from ours. One can check
that the volume of $X$ in [4] corresponds to $\frac{n!}{(4\pi)^n}
\v(X)$ in our notation.

One consequence of Hirzebruch proportionality is a formula for the
dimension of the space $V_{\ell}$ of cusp forms of weight $\ell$.
 By definition, $V_{\ell}$ is the space of
sections of $L^{\otimes \ell}$ which vanish on $E$. In other
words,
$$V_{\ell} := H^0( \bar{X}, {\cal O}(\ell L -E)).$$
Mumford showed that the formula for the dimension of $V_{\ell}$ in
the case of compact $X$ continues to hold for non-compact $X$ with
an error term of degree bounded by the dimension of $X^*\setminus
X$. More precisely,

\medskip
{\bf Proposition 2 [4, Corollary 3.5]} {\it Let $$P(\ell) :=
h^0({\bf P}_n, {\cal O}(\ell(n+1))) = \frac{(n \ell + n +
\ell)!}{n! (n \ell + \ell)!}.$$ Then there exists a constant $P_0$
such that for all $\ell \geq 2$,
$$ \dim V_{\ell} = \frac{n!}{(4 \pi)^n}\v(X)  P(\ell -1) + P_0.$$}

\medskip
An immediate consequence is

\medskip
{\bf Corollary 1} {\it For any $\ell \geq 2$, $\dim V_{\ell +1}
 > \dim V_{\ell}.$ In particular, $V_3 \neq 0$. }

\section{Proof of Theorem 2}

\medskip
To prove Theorem 2, we need the following two lemmas.

\medskip
{\bf Lemma 1} {\it Recall that $E_1, \ldots, E_k$ are the
components of $E = \bar{X} \setminus X$. For each $1 \leq i \leq
k$, there exists $\sigma_i \in H^0(\bar{X}, {\cal O}(2 L))$ such
that
$$\sigma_i |_{E_i} \neq 0, \mbox{ but } \sigma_i |_{E_j} = 0 \mbox{ for each }
j \neq i.$$}

\medskip
{\it Proof of Lemma 1}. Consider the short exact sequence on
$\bar{X}$,
$$ 0 \longrightarrow {\cal O}(2L-E) \longrightarrow
{\cal O}(2L) \longrightarrow {\cal O}(2L)|_{E} \longrightarrow
0.$$ Since $L= K_{\bar{X}} +E$ is nef and big, Kawamata-Viehweg
vanishing gives $$H^1(\bar{X}, {\cal O}(2L-E)) = H^1(\bar{X},
{\cal O}(K_{\bar{X}} + L)) = 0.$$ Thus we have the surjectivity of
the restriction map
$$H^0(\bar{X}, {\cal O}(2L)) \rightarrow H^0(E, {\cal
O}(2L)|_{E}).$$ Since $E_i$ is contracted by $\psi: \bar{X}
\rightarrow X^*$, the line bundle $L|_{E_i}$ is trivial. So we
have the surjectivity of
$$H^0(\bar{X}, {\cal O}(2L)) \rightarrow \bigoplus_{i=1}^k
H^0(E_i, {\cal O}_{E_i}) $$ from which Lemma 1 follows. $\Box$

\medskip
{\bf Lemma 2} {\it Suppose $V_{\ell} \neq 0$. Then $\dim V_{\ell +
2} - \dim V_{\ell} \geq k-1$. }

\medskip
{\it Proof of Lemma 2}. Recall that elements of $ V_{\ell}$ are
sections of $L^{\otimes \ell}$ which vanish on $E$. Choose $v \in
V_{\ell}$ such that the vanishing order of $v$ along $E_1$ is the
highest among all non-zero elements of $V_{\ell}$. Fix a basis $\{
v_1, \ldots, v_m\}$ of $V_{\ell}$ with $m = \dim V_{\ell}$.
Consider the following $(m+ k-1)$ elements of $V_{\ell +2}$.
$$\sigma_2 \cdot v, \ldots, \sigma_k \cdot v, \sigma_1 \cdot v_1,
\ldots, \sigma_1 \cdot v_m$$ where $\sigma_1, \ldots, \sigma_k$
are as in Lemma 1. We claim that they are linearly independent.
Suppose
$$ \sum_{j=2}^k a_j (\sigma_j \cdot v) + \sum_{i=1}^m b_i
(\sigma_1 \cdot v_i) =0$$ for some complex numbers $a_j, b_i$.
Then $$ (\sum_{j=2}^k a_j \sigma_j) \cdot v = - \sigma_1 \cdot w$$
for $w= \sum_{i=1}^m b_i v_i \in V_{\ell}$. The left hand side has
vanishing order along $E_1$ strictly higher than that of $v$.
Since the vanishing order of non-zero $w$ along $E_1$ can't be
bigger than that of $v$, we see that $w=0$. This yields $a_j = b_i
=0$ for all $2 \leq j \leq k$ and $1 \leq i \leq m$. This proves
the claim. Lemma 2 follows immediately from the claim. $\Box$

\medskip
{\it Proof of Theorem 2}. {}From Corollary 1 and Lemma 2, we see
that
$$\dim V_5 -\dim V_3 \geq k-1.$$ By Proposition 2,
$$\dim V_5 - \dim V_3 = \frac{n!}{(4 \pi)^n}\v(X) (P(4) - P(2))>0 .$$
Thus  $$\v(X) \geq \frac{(4 \pi)^n}{n! (P(4) - P(2))} (k-1).$$

\smallskip
 As quoted in [3, p.179], Gromov's generalization of
Gauss-Bonnet says
$$  \v(X) = \frac{(-4\pi)^n}{(n+1)!} e(X)$$ where $e(X)$ denotes the topological
Euler number of $X$. This implies
 $$ \v(X) \geq \frac{(4 \pi)^n}{(n+1)!}.$$ Thus when $k \leq
\frac{P(4) - P(2)}{n+1}$,
 $$
\v(X) \geq \frac{(4 \pi)^n}{n! (P(4) - P(2))}  k$$ and the
statement of  Theorem 2 holds automatically.

\smallskip
When $k \geq \frac{P(4) - P(2)}{n+1}$,$$ k-1 \geq (1-
\frac{n+1}{P(4) - P(2)}) k.$$ Thus
$$\v(X) \geq \frac{(4 \pi)^n}{n! (P(4) - P(2))} (k-1) \geq
\frac{(4 \pi)^n}{ n! (P(4) - P(2))} (1- \frac{n+1}{P(4) - P(2)})
k$$  which proves the theorem. $\Box$.

\bigskip

 {\bf References}

\medskip
[1] A. Ash, D. Mumford, M. Rapoport and Y. Tai, {\it Smooth
compactification of locally symmetric varieties,} Math. Sci.
Press, 1975

[2] W. L. Baily and A. Borel, Compactification of arithmetic
quotients of bounded symmetric domains, {\it Annals of Math.} {\bf
84} (1966), 442-528

 [3] S. Hersonsky and F. Paulin, On the volumes of complex
 hyperbolic manifolds, {\it Duke Math. J.} {\bf 84} (1996), 719-737

[4] D. Mumford, Hirzebruch's proportionality theorem in the
non-compact case, {\it Invent. math.} {\bf 42} (1977), 239-272

 [5] J. R. Parker, On the volumes of cusped, complex
hyperbolic manifolds and orbifolds, {\it Duke Math. J.} {\bf 94}
(1998), 433-464

[6] Y.-T. Siu and S.-T. Yau,  Compactifications of negatively
curved complete K\"ahler manifolds of finite volume, {\it  Ann.
Math. Stud.} {\bf 102} (1980), 363-380

 [7] W.-K. To,  Total geodesy of proper holomorphic immersions
 between complex hyperbolic space forms of finite volume, {\it Math. Ann.} {\bf
 297} (1993), 59-84

\bigskip
Korea Institute for Advanced Study

207-43 Cheongryangri-dong

Seoul, 130-722,  Korea

jmhwang@kias.re.kr

\end{document}